\documentclass[dvips,12pt,a4paper]{article}
\usepackage{amsmath,amssymb}        % moduli per le estensioni ed i fonts AMS
\usepackage{latexsym}
\usepackage{graphicx}
\usepackage{psfrag}
\usepackage{url}

\newtheorem{theorem}{Teorema}
\def\RR{\hbox{I\kern-.2em\hbox{R}}}
\def\ds{\displaystyle}

\def\tw{.85\textwidth}

\title{The Non-Iterative Transformation Method}

\author{Riccardo Fazio  \\
Department of Mathematics, Computer Science,\\ Physical Sciences and Earth Sciences,\\
University of Messina,\\
Viale F. Stagno D'Alcontres 31, 98166 Messina, Italy. \\ e-mail: rfazio@unime.it \\ home-page: http://mat521.unime.it/fazio.}
\begin{document}

\maketitle

%\contributed{$^\dagger$ These authors contributed equally to this work.}

% Contact information of the corresponding author (Add [2] after \corres if there are more than one corresponding author.)
%\corres{}

% Abstract (Do not use inserted blank lines, i.e. \\)
\abstract{The Blasius flow is the idealized flow of a viscous fluid past an infinitesimally thick, semi-infinite flat plate. 
The definition of a non-iterative transformation method for the celebrated Blasius problem is due to T{\"o}pfer and dates more than a century ago.
Here we define a non-iterative transformation method for Blasius equation with a moving wall, a slip flow condition or a surface gasification.
The defined method allows us to deal with classes of problems in boundary layer theory that, depending on a parameter, admit multiple or no solutions.
This approach is particularly convenient when the main interest is on the behaviour of the considered models with respect to the involved parameter.
The obtained numerical results are found to be in good agreement with those available in literature.    
}

% Keywords: add 3 to 10 keywords

% The fields PACS, MSC, and JEL may be left empty or commented out if not applicable
%\PACS{}
%\MSC{}
%\JEL{}

% If this is an expanded version of a conference paper, please cite it here: enter the full citation of your conference paper, and add $^\dagger$ in the end of the title of this article.
%\conference{}

%%%%%%%%%%%%%%%%%%%%%%%%%%%%%%%%%%%%%%%%%%
% For journal Data:

%\dataset{DOI number or link to the deposited data set in cases where the data set is published or set to be published separately. If the data set is submitted and will be published as a supplement to this paper in the journal Data, this field will be filled by the editors of the journal. In this case, please make sure to submit the data set as a supplement when entering your manuscript into our manuscript editorial system.}
%\datasetlicense{license under which the data set is made available (CC0, CC-BY, CC-BY-SA, CC-BY-NC, etc.)}

%%%%%%%%%%%%%%%%%%%%%%%%%%%%%%%%%%%%%%%%%%

\section{Introduction}\label{S:History}

At the beginning of the last century Prandtl \cite{Prandtl:1904:UFK} put the foundations of
boundary-layer theory providing the basis for the unification of two, at that time seemingly incompatible, sciences: namely, theoretical hydrodynamics and hydraulics.
Boundary-layer theory has found its main application in calculating the skin-friction drag which acts on a body as it is moved through a fluid: for example the drag of an airplane wing, of a turbine blade, or a complete ship \cite{Schlichting:1979:BLT}.
The Blasius flow is the idealized flow of a viscous fluid past an infinitesimally thick, semi-infinite
flat plate.
Recently, Boyd \cite{Boyd:2008:BFC} uses Blasius problem as an example were some good analysis, before the computer invention, allowed researchers of the past to solve problems, governed by partial differential equations, that might be otherwise impossible to face.

Blasius problem \cite{Blasius:1908:GFK} is the simplest nonlinear boundary layer problem.
A study by Boyd point out how this particular problem has arisen the interest of prominent scientist, like H. Weyl, J. von Neumann, M. Van Dyke, etc., see Table 1 in \cite{Boyd:1999:BFC}. 
The main reason for this interest is due to the hope that any approach developed for this epitome can be extended to more difficult hydrodynamics problems.

Blasius main interest was to compute, without worrying about existence or uniqueness of its boundary value problem (BVP) solution -- it was Weyl who proved in \cite{Weyl:1942:DES} that Blasius problem has one and only one solution --, the value of  $\lambda$ the so-called shear stress.
To compute this value, Blasius used a formal series solution around $\eta=0$ and an asymptotic expansions for large values of $\eta$, adjusting the constant $\lambda$ so as to connect both expansions in a middle region.
In this way, Blasius obtained the (erroneous) bounds $ 0.3315 < \lambda < 0.33175$.

A few years later, T{\"o}pfer \cite{Topfer:1912:BAB} revised the work by Blasius and solved numerically the Blasius equation with suitable initial conditions and the classical order-four Runge-Kutta method.
He then arrived, without detailing his computations, at the value
$\lambda \approx 0.33206$, contradicting the bounds reported by Blasius.  

Thereafter, the quest for a good approximation of $\lambda$ was a main concern.
This is seldom the case for the most important problems of applied mathematics:
at the first study everyone would like to know if there is a method to solve a given problem, but, as soon as a problem is solved, then we would like to know how accurate is the computed solution and whether there are different methods that can provide a solution with less effort. 

By using a power series, Bairstow \cite{Bairstow:1925:SF} reports $\lambda \approx 0.335$, and Goldstein \cite{Goldstein:1930:CSS} obtains $\lambda \approx 0.332$ or, using a finite difference method, Falkner \cite{Falkner:1936:MNS} finds $\lambda \approx 0.3325765$, and Howarth \cite{Horwarth:1938:SLB} yields $\lambda \approx 0.332057$.
Fazio \cite{Fazio:1992:BPF}, using a free boundary formulation of the Blasius problem, finds $\lambda \approx 0.332057336215$.
Boyd \cite{Boyd:1999:BFC} uses T{\"o}pfer's algorithm to obtain the accurate value $\lambda \approx 0.33205733621519630$.
By the Adomain's decomposition method Abbasbandy \cite{Abbasbandy:2007:NSB} finds $\lambda \approx 0.333329$, whereas a variational iteration method with Pad\'e approximants allows Wazwaz \cite{Wazwaz:2007:VIM} to calculate, the imprecise value, $\lambda \approx 0.3732905625$.
Tajvidi et al. \cite{Tajvidi:2008:MRL} apply modified rational Legendre functions to get a value of  $\lambda \approx 0.33209$.

At the turning of this new century, as the number of applications of microelectronics devices increases, boundary-layer theory has found a renewal of interest within the study of gas and liquid flows at the micro-scale regime, see, for instance, Gad el Hak \cite{Gad-el-Hak:1999:FMM} or Martin and Boyd \cite{Martin:2001:BBL}.

\section{Fluid flow on a flat plate}\label{S:Fluid}

\begin{figure}[!hbt]
	\centering
\psfrag{y}[][]{$ y $} 
\psfrag{z}[][]{$ z $} 
\psfrag{O}[][]{$ O $}
\psfrag{V}[][]{$ V_\infty $} 
\includegraphics[width=\tw]{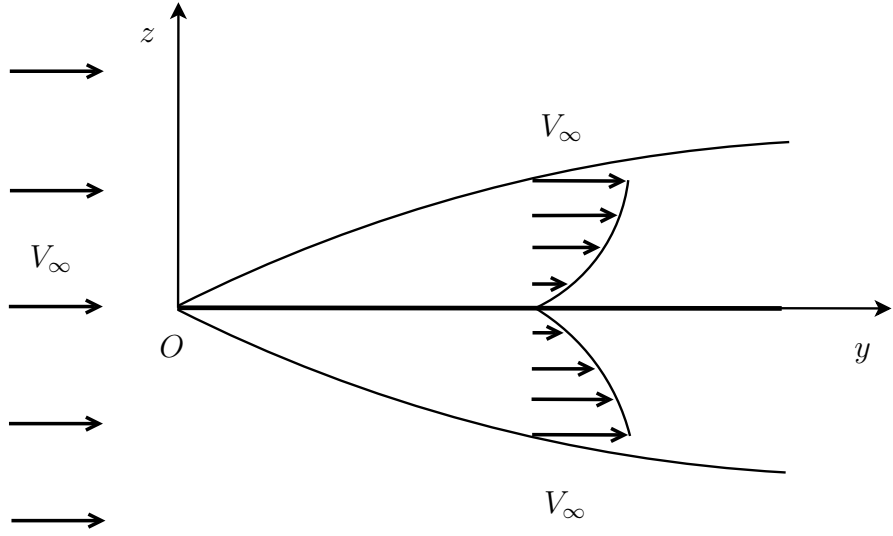}
\caption{Sketch of boundary layer over a thin plate.}
	\label{fig:plate}
\end{figure}

The model describing
the steady plane flow of a fluid past a thin plate\index{fluid!on a thin plate}%
, provided the boundary layer assumptions are verified (the flow has a very thin layer attached to the plate and $ v \gg w $), is given by
\begin{align}\label{PDE-model}
& {\displaystyle \frac{\partial v}{\partial y}} +
{\displaystyle \frac{\partial w}{\partial z}} = 0  \nonumber \\
& v {\displaystyle \frac{\partial v}{\partial y}} +
w {\displaystyle \frac{\partial v}{\partial z}} = \nu
{\displaystyle \frac{\partial^2 v}{\partial z^2}}  \\
& v(y, 0) = w(y, 0) = 0 \ , \quad v(y, z) \rightarrow V_{\infty}
\quad \mbox{as} \quad z \rightarrow \infty \ , \nonumber 
\end{align}
where the governing differential equations, namely conservation of mass and momentum, are the 
steady-state 2D Navier-Stokes
equations under the boundary layer approximations,
$ v $ and $ w $ are the velocity components of the fluid in the
$ y $ and $ z $ direction, $ V_{\infty} $
represents the main-stream velocity, see the draft in figure \ref{fig:plate}, and $ \nu $ is the viscosity of the fluid.
The boundary conditions at $ z = 0 $ are based on the
assumption that neither slip nor mass
transfer are permitted at the plate whereas the remaining boundary condition
means that the velocity $ v $ tends to the main-stream velocity
$ V_{\infty} $ asymptotically.

In order to study this problem, it is convenient to introduce a
potential (stream function) $ \psi(y, z) $ defined by
\begin{equation}
v = {\displaystyle \frac{\partial \psi}{\partial z}} \ , \qquad
w = - {\displaystyle \frac{\partial \psi}{\partial y}} \ .
\end{equation}
The physical motivation for introducing this function is that constant
$ \psi $ lines are steam-lines.
The mathematical motivation for introducing such a new variable is that the equation of continuity is satisfied identically, and we have to deal only with the transformed momentum equation.
In fact, introducing the stream function the problem can be rewritten as follows
\begin{align}\label{stream-model}
& \nu {\displaystyle \frac{\partial^3 \psi}{\partial z^3}} +
{\displaystyle \frac{\partial \psi}{\partial y}}
{\displaystyle \frac{\partial^2 \psi}{\partial z^2}} -
{\displaystyle \frac{\partial \psi}{\partial z}}
{\displaystyle \frac{\partial^2 \psi}{\partial y \partial z}}
 = 0  \nonumber \\
& {\displaystyle \frac{\partial \psi}{\partial y}}(y, 0) =
{\displaystyle \frac{\partial \psi}{\partial z}} (y, 0) = 0  \\
& {\displaystyle \frac{\partial \psi}{\partial z}} (y, z)
\rightarrow V_{\infty}  \quad \mbox{as} \quad
z \rightarrow \infty \ . \nonumber
\end{align}

\subsection{Blasius problem}

Blasius \cite{Blasius:1908:GFK} used the following similarity transformation
\begin{equation}
 \eta = z \left(\frac{V_\infty}{\nu y} \right)^{1/2} \ , \qquad  f(\eta) = \psi(y, z)
\left(\nu y V_\infty \right)^{-1/2} \ , 
\end{equation}
that reduces the partial differential model (\ref{stream-model}) to \index{Blasius problem}%
\begin{align}\label{eq:Blasius}
&{\displaystyle \frac{d^3 f}{d \eta^3}} + \frac{1}{2} f
{\displaystyle \frac{d^{2}f}{d\eta^2}} = 0 \nonumber \\[-1.5ex]
& \\[-1.5ex]
&f(0) = {\displaystyle \frac{df}{d\eta}}(0) = 0 \ , \quad
{\displaystyle \frac{df}{d\eta}}(\eta) \rightarrow 1 \quad \mbox{as}
\quad \eta \rightarrow \infty \ , \nonumber
\end{align}
i.e., a boundary value problem (BVP) defined on a semi-infinite interval.
Blasius solved this BVP by patching a power series to an asymptotic approximation at some finite value of $ \eta $.

\section{T{\"o}pfer transformation}\label{Toepfer}

In order to clarify T{\"o}pfer \cite{Topfer:1912:BAB} derivation of a further transformation of variables that reduces the BVP into an initial value problem (IVP) we consider the derivation of the series expansion solution. %see for instance Na \cite[pp. 140-142]{Na:1979:CME}.
Of course, some of the coefficients of the series can be evaluated by imposing the boundary conditions at $ \eta = 0 $.
Moreover, we set
\begin{equation}
    \lambda = \frac{d^2f}{d\eta^2}(0) 
\end{equation} 
where $ \lambda $ is a constant different from zero.
So that, we look for a series solution defined as
\begin{equation}
    f(\eta) = \frac{\lambda}{2} \eta^2 + \sum_{n=3}^\infty C_n \eta^n
\end{equation}
where the coefficients $ \lambda $ and $ C_n $, for $ n = 3, 4, \dots $, are constants to be determined.
In fact, the boundary values at the plate surface, at $\eta = 0$,
require that $ C_0 = C_1 = 0 $, and we also have $ C_2 = \lambda/2 $ by the definition of $ \lambda $.
Now, we substitute this series expansion into the governing differential equation, whereupon we find 
\begin{multline} 
\sum_{n=3}^\infty n (n-1)(n-2)C_n \eta^{n-3} \\
 + \frac{1}{2} 
\left(\frac{\lambda}{2} \eta^2 + \sum_{n=3}^\infty C_n \eta^n\right)
\left[\lambda + \sum_{n=3}^\infty n (n-1) C_n \eta^{n-2} \right] = 0 
\end{multline}
or in expanded form 
\begin{multline}
 \left[3 \; 2 \; C_3 \right] 
+ \left[4 \; 3 \; 2 \; C_4 \right] \eta 
+ \left[5 \; 4 \; 3 \; C_5 + \frac{1}{2} \; 2 \; \frac{\lambda}{2} \frac{\lambda}{2}\right] \eta^2 \\
 + \left[6 \; 5 \; 4 \; C_6 + \frac{1}{2} \; 2 \; \frac{\lambda}{2} C_3 
 + \frac{1}{2} \; \frac{\lambda}{2} \; 3 \; 2 \; C_3\right] \eta^3 \\
 + \cdots = 0 \ .
\end{multline}
According to a standard approach, we have to require that all coefficients of the powers of $ \eta $ to be zero.
It is an easy matter to compute the coefficients of the series expansion in terms of $ \lambda $:
\begin{align*}
C_3 &= C_4 = 0 \ , \qquad C_5 = -\frac{\lambda^2}{2\; 5!} \nonumber \\
C_6 &= C_7 = 0 \ , \qquad C_8 = 11\frac{\lambda^3}{2^2\; 8!} \nonumber \\
C_9 &= C_{10} = 0 \ , \qquad C_{11} = -375\frac{\lambda^4}{2^3\; 11!} \nonumber \\
&\mbox{and so on ...} \nonumber 
\end{align*}
The solution can be written as
\begin{equation}
f = \frac{\lambda \eta^2}{2} - \frac{\lambda^2 \eta^5}{2\; 5!} 
+ \frac{11 \; \lambda^3 \eta^8}{2^2\; 8!} 
- \frac{375 \; \lambda^4 \eta^{11}}{2^3\; 11!} + \cdots 
\end{equation}
where the only unknown constant is $ \lambda $.
In principle, $ \lambda $ can be determined by imposing the boundary condition at the second point, but, in this case, this cannot be done because the left boundary condition is given at infinity.
However, by modifying the powers of $ \lambda $ we can rewrite the series expansion as
\begin{equation} 
\lambda^{-1/3} f = \frac{\left(\lambda^{1/3} \eta\right)^2}{2} - \frac{\left(\lambda^{1/3} \eta\right)^5}{2\; 5!}  
+ \frac{11 \; \left(\lambda^{1/3} \eta\right)^8}{2^2\; 8!} 
- \frac{375 \; \left(\lambda^{1/3} \eta\right)^{11}}{2^3\; 11!} + \cdots 
\end{equation} 
which suggests a transformation of the form \index{scaling transformation}%
\begin{equation}\label{eq:scalinv:Blasius}
f^* = \lambda^{-1/3} f \ , \qquad \eta^* = \lambda^{1/3} \eta  \ .   
\end{equation}
In the new variables the series expansion becomes 
\begin{equation}
f^* = \frac{\eta^{*2}}{2} - \frac{\eta^{*5}}{2\; 5!} 
+ \frac{11 \; \eta^{*8}}{2^2\; 8!} 
- \frac{375 \; \eta^{*11}}{2^3\; 11!} + \cdots 
\end{equation}
which does not depend on $ \lambda $.
We notice that the governing differential equation and the initial conditions at the free surface, at $\eta = 0$, are left invariant by the new variables defined above.
Moreover, the first and second order derivatives transform in the following way
\begin{equation}
\frac{d f^*}{d \eta^{*}} = \lambda^{-2/3} \frac{d f}{d \eta} \ ,   
\qquad
\frac{d^2 f^*}{d \eta^{*2}} = \lambda^{-1} \frac{d^2 f}{d \eta^{2}} \ .  
\end{equation}
As a consequence of the definition of $ \lambda $ we have
\begin{equation}
\frac{d^2 f^*}{d \eta^{*2}} (0) = 1 \ ,   
\end{equation} 
and this explains why in these variables the series expansion does not depend on $ \lambda $.
Furthermore, the value of $ \lambda $ can be found on condition that we have an approximation for $ \frac{d f^*}{d \eta^{*}}(\infty) $, say $ \frac{d f^*}{d \eta^{*}}(\eta_{\infty}) $ where $\eta^*_{\infty}$ is a suitable truncated boundary.
In fact, by the above relation we get
\begin{equation}\label{eq:lambda:Blasius}
\lambda = \left[ \frac{d f^*}{d \eta^{*}}(\eta^*_{\infty}) \right]^{-3/2} \ .   
\end{equation} 

From a numerical viewpoint, BVPs must be solved within the computational domain simultaneously (a \lq \lq stationary\rq \rq \ problem), whereas IVPs can be solved by a stepwise procedure (an \lq \lq evolution\rq \rq \ problem).
Somehow, numerically, IVPs are easier than BVPs.  

\subsection{T{\"o}pfer algorithm}

\begin{figure}[!hbt]
	\centering
\psfrag{e}[][]{$ \eta^*, \eta $} 
\psfrag{f2star}[][]{$ \frac{df^*}{d\eta^*} $} 
\psfrag{f2}[][]{$ \frac{df}{d\eta} $} 
\psfrag{f3star}[][]{$ \frac{d^2f^*}{d\eta^{*2}} $} 
\psfrag{f3}[][]{$ \frac{d^2f}{d\eta^2} $} 
\includegraphics[width=\tw]{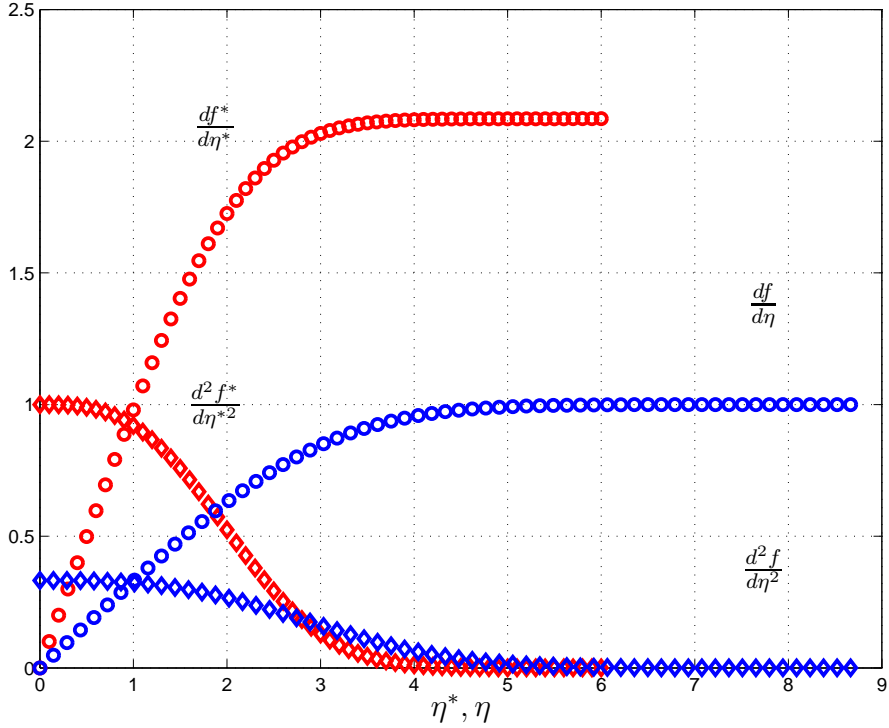}
\caption{Blasius solution by a non-ITM for $\eta^* \in [0,6]$ we get $\lambda_2 = 0.332057687$.} 
	\label{fig:Blasius2}
\end{figure}

Let us list the steps necessary to solve the Blasius problem by the T{\"o}pfer algorithm.
In this way, we define a non-iterative (I)TM.
We have to:\index{T{\"o}pfer algorithm}%
\begin{enumerate}
    \item solve the auxiliary IVP \index{IVP!auxiliary}%
\begin{align}\label{eq:Blasius2}
&{\displaystyle \frac{d^3 f*}{d \eta^{*3}}} + \frac{1}{2} f^*
{\displaystyle \frac{d^{2}f^*}{d\eta^{*2}}} = 0 \nonumber \\[-1.5ex]
& \\[-1.5ex]
&f^*(0) = {\displaystyle \frac{df^*}{d\eta^*}}(0) = 0, \qquad
{\displaystyle \frac{d^2f^*}{d\eta^{*2}}}(0) = 1  \nonumber
\end{align}
and, in particular, get an approximation for
$ \frac{d f^*}{d \eta^{*}}(\infty) $;
\item
compute $ \lambda $ by equation (\ref{eq:lambda:Blasius});
\item
obtain $ f(\eta) $, ${\displaystyle \frac{df}{d\eta}}(\eta)$, and ${\displaystyle \frac{d^2f}{d\eta^2}}(\eta)$ by the inverse transformation of (\ref{eq:scalinv:Blasius}).
\end{enumerate}
Indeed, T{\"o}pfer solved the IVP for the Blasius equation once.
At large but finite $ \eta_j^* $, ordered so that $ \eta_j^* < \eta_{j+1}^* $, we can compute by equation (\ref{eq:lambda:Blasius}) the corresponding $ \lambda_j $.  
If two subsequent values of $ \lambda_j $ agree within a specified accuracy, then $ \lambda $ is approximately equal to the common value of the $ \lambda_j $, otherwise, we can march to a larger value of $ \eta $ and try again.
Using the classical fourth order Runge-Kutta \index{Runge-Kutta method}%
 method and a grid step $ \Delta \eta^* = 0.1$
T{\"o}pfer was able to determine $ \lambda $ with an error less than $ 10^{-5} $. 
He used the two truncated boundaries $\eta_1^* = 4$ and $\eta_2^* = 6$.
We notice that for $\eta^* \in [0,4]$ we get $\lambda_1 = 0.333233336$.
% see the figures \ref{fig:Blasius1}-\ref{fig:Blasius2} were we have reproduced T{\"o}pfer computations.
%\begin{figure}[!hbt]
%	\centering
%\psfrag{e}[][]{$ \eta^*, \eta $} 
%\psfrag{f2star}[][]{$ \frac{df^*}{d\eta^*} $} 
%\psfrag{f2}[][]{$ \frac{df}{d\eta} $} 
%\psfrag{f3star}[][]{$ \frac{d^2f^*}{d\eta^{*2}} $} 
%\psfrag{f3}[][]{$ \frac{d^2f}{d\eta^2} $} 
%\includegraphics[width=\tw]{Figs/TopferA} 
%\caption{Blasius solution by a non-ITM for $\eta^* \in [0,4]$ we get $\lambda_1 = 0.333233336$} 
%	\label{fig:Blasius1}
%\end{figure}
%In the figures \ref{fig:Blasius1}-\ref{fig:Blasius2}
In figure \ref{fig:Blasius2} we plot the more accurate numerical solution obtained by T{\"o}pfer's algorithm defined above.
We notice that this figure shows the solutions of the auxiliary IVP (\ref{eq:Blasius2}) and of the BVP (\ref{eq:Blasius}).

\section{Rubel error analysis}

The boundary condition at infinity is certainly not suitable for a numerical treatment.
This condition has usually been replaced by the same condition applied at a truncated boundary, see Collatz \cite[pp 150-151]{Collatz} or Fox \cite[p. 92]{Fox}.
In the truncated boundary formulation $ f_M(\eta) $ is defined by
\begin{eqnarray}
& {\displaystyle \frac{d^3 f_M}{d \eta^3}} + f_M
{\displaystyle \frac{d^{2}f_M}{d\eta^2}} = 0 \nonumber \\[-1.5ex]
\label{p7} \\[-1.5ex]
& f_M(0) = {\displaystyle \frac{df_M}{d\eta}}(0) = 0, \qquad
{\displaystyle \frac{df_M}{d\eta}}(M) = 1  \nonumber
\end{eqnarray}
where $ M $ represents the truncated boundary. 
It is evident that also in (\ref{p7}) the governing DE and the two boundary
conditions at the origin are left invariant by the scaling transformation (\ref{eq:scalinv:Blasius}).

For the error related to the truncated boundary solution $ f_M(\eta) $ defined by
\begin{equation}
e(\eta) = |f(\eta) - f_M(\eta)|, \qquad \eta \in [0, M] \ ,
\end{equation}
the following theorem holds true.

\begin{theorem} {\bf \ (due to Rubel \cite{Rubel:1955:EET}).} A truncated boundary formulation of the Blasius problem introduces an error which verifies the following inequality
\begin{equation}
e(\eta) \leq M {\displaystyle \frac{d^2 f_M}{d\eta^2}}(M)
[f_M(M)]^{-1}  \quad .
\end{equation}
\end{theorem}

\medskip
\noindent
{\bf Outline of the proof.} 
As proved by Weyl \cite{Weyl:1942:DES}, it is true that
\begin{equation}\label{Weylin}
{\displaystyle \frac{d^2 f}{d\eta^2}}(\eta) > 0, \qquad
\mbox{for \ \ } \eta \in (0, \infty)  \ .
\end{equation}
By (\ref{Weylin}) and taking into account the boundary conditions in (\ref{eq:Blasius}), we have that
\begin{equation}
{\displaystyle \frac{df}{d\eta}}(\eta) \ \ \mbox{and} \ \ f(\eta) \ \
\mbox{are increasing functions on} \ \ \eta \in (0, \infty) \quad .
\end{equation}
As a consequence the function $ \lambda^2 {\displaystyle \frac{df}{d\eta}}(\lambda M) $ is zero for $ \lambda = 0 $, increases with $ \lambda $ and tends to infinity as $ \lambda \rightarrow \infty $. 
For some value of $ \lambda \in (0, \infty) $ we must have
\begin{equation}
\lambda^2 {\displaystyle \frac{df}{d\eta}}(\lambda M) = 1  \quad .
\end{equation}
This value verifies $ \lambda > 1 $ because $ \lambda M $ is a
finite value, $ {\displaystyle \frac{df}{d\eta}}(\eta) $ is an
increasing function and
$ {\displaystyle \frac{df}{d\eta}}(\eta) \rightarrow 1 \ \ \mbox{as}
\ \ \eta \rightarrow \infty $.
For this particular value of $ \lambda $, due to the scaling properties, we have found that
\begin{equation}
f_M(\eta) = \lambda f(\lambda \eta)
\end{equation}
because $ \lambda f(\lambda \eta) $ verifies the BVP (\ref{p7}) that defines $ f_M(\eta) $ uniquely.

Hence, the error for $ \eta \in [0, M] $ is given by
\begin{equation}
e(\eta) = |\lambda f(\lambda \eta) - f(\eta)| \leq
|(\lambda -1) f(\lambda \eta)| + |f(\lambda \eta) - f(\eta)| \quad .
\end{equation}
By applying the mean value theorem of differential calculus and taking
into account that $ {\displaystyle \frac{df}{d\eta}}(\eta) \leq 1 $
we get the relations $ f(\lambda \eta) \leq \lambda \eta $ and
$ |f(\lambda \eta) - f(\eta)| \leq (\lambda - 1) \eta $. As a result
\begin{equation}
e(\eta) \leq M (\lambda^2 -1), \qquad \eta \in [0, M]
\end{equation}
where $ \lambda^2 -1 > 0 $ because $ \lambda > 1 $.
Naturally,
$ {\displaystyle \frac{df_M}{d\eta}}(\eta \rightarrow \infty) =
\lambda^2 $, so that
\begin{eqnarray*}
& \begin{array}{ll}
\lambda^2 -1 &= {\displaystyle \frac{df_M}{d\eta}}(\eta \rightarrow \infty)
- {\displaystyle \frac{df_M}{d\eta}}(M) \\
& = {\displaystyle \int_{M}^{\infty}}
{\displaystyle \frac{d^2 f_M}
{d\eta^2}} d\eta \quad .
\end{array}
\end{eqnarray*}
To complete the proof Rubel used some manipulations,
involving a first integral of the governing differential equation, to find that
\begin{equation}
\lambda^2 -1 \leq {\displaystyle \frac{d^2 f_M}{d\eta^2}}(M)
[f_M(M)]^{-1} \ .
\end{equation}
\hfill $ \Box $

\noindent
{\bf Remark.} As a consequence of this theorem in order to control the error we can modify either the value of $ M $ or the value of $ {\displaystyle \frac{d^2 f_M}{d\eta^2}}(M) $. 
Classically the value of $ M $ has been chosen to this end.
The above Theorem shows that the error is directly proportional to $ M $. 
In this context Fazio defines a free boundary formulation of the Blasius problem where the second order derivative of the solution with respect to $ \eta $ at the free boundary can be chosen as small as possible, see \cite{Fazio:1992:BPF} for details.
The free boundary can be interpreted as a truncated boundary, see Fazio \cite{Fazio:2002:SFB}.

%%%%%%%%%%%%%%%%%%%%%%%%%%%%%%%%%%%%%%%%%%

\section{Similarity analysis}

Motivated by several problems in boundary layer theory, let us consider the class of BVPs defined by
\begin{align}\label{eq:class1}
&{\displaystyle \frac{d^3 f}{d \eta^3}} = \phi\left(\eta, f,
{\displaystyle \frac{df}{d\eta}, \frac{d^{2}f}{d\eta^2}}\right) \nonumber \\[-1.5ex]
& \\[-1.5ex]
& f(0) = a \, \quad {\displaystyle \frac{df}{d\eta}}(0) = b + c {\displaystyle \frac{d^2f}{d\eta^2}}(0) \ , \quad
{\displaystyle \frac{df}{d\eta}}(\eta) \rightarrow d \quad \mbox{as}
\quad \eta \rightarrow \infty \ , \nonumber
\end{align}
where $a$, $b$, $c$ and $d$ are given constants, with $d \ne 0$.
Introducing  the scaling group
\begin{equation}\label{eq:scaling1}
f^* = \lambda f \ , \qquad \eta^* = \lambda^{\delta} \eta \ ,   
\end{equation}
we require the invariance of (\ref{eq:class1}), but the asymptotic boundary condition so that $\delta \ne 1$, with respect to (\ref{eq:scaling1}).
The requested invariance is granted on condition that $a=b=c=0$ and
\begin{equation}\label{eq:condition}
\phi\left(\eta, f,{\displaystyle \frac{df}{d\eta}, \frac{d^{2}f}{d\eta^2}}\right) =
\eta^{1-3\delta} \Phi\left(\eta^{1/\delta} f, \eta^{(1-\delta)/\delta}{\displaystyle \frac{df}{d\eta}}, \eta^{(1-2\delta)/\delta}{\displaystyle \frac{d^{2}f}{d\eta^2}}\right)  \ .   
\end{equation}
As a consequence of the above scaling invariance we can define a non-ITM.

\subsection{The non-iterative algorithm}

In order to define the numerical method for the characterized class of problems we have to consider the auxiliary IVP
\begin{align}\label{eq:IVP1}
&{\displaystyle \frac{d^3 f^*}{d \eta^{*3}}} = \eta^{*(1-3\delta)} \Phi\left(\eta^{*(1/\delta)} f^*, \eta^{*(1-\delta)/\delta}{\displaystyle \frac{df^*}{d\eta^*}}, \eta^{*(1-2\delta)/\delta}{\displaystyle \frac{d^{2}f^*}{d\eta^{*2}}}\right)  \nonumber \\[-1.5ex]
& \\[-1.5ex]
& f^*(0) = {\displaystyle \frac{df^*}{d\eta^*}}(0) = 0 \ , \qquad {\displaystyle \frac{d^2f^*}{d\eta^{*2}}}(0) = p \ .  \nonumber
\end{align}
where $p$ is defined by the user, we usually set $p=\pm 1$, but it is also possible to consider different values.
For instance, in Fazio \cite{Fazio:1992:BPF} for the Blasius problem we used $ p = 1000$.
We have to solve (\ref{eq:IVP1}) on $[0, \eta_{\infty}^*]$, where $\eta_{\infty}^*$ is a suitable truncated boundary chosen under the condition
\begin{equation}\label{eq:asymcond}
{\displaystyle \frac{df^*}{d\eta^*}}(\eta_{\infty}^*) \approx {\displaystyle \frac{df^*}{d\eta^*}}(\infty) \ .
\end{equation}
As $d \ne 0$, we have
\begin{equation}\label{eq:lambda1}
\lambda = \left[ {\ds \frac{\frac{d f^*}{d \eta^{*}}(\eta_{\infty}^*)}{d}} \right]^{1/(1-\delta)} \ .   
\end{equation} 
Computed the value of $\lambda$ we can apply the inverse transformation of (\ref{eq:scaling1}) to get
\begin{align}\label{eq:rescale1} 
& \eta = \lambda^{-\delta} \eta^* \ , \quad f(\eta) = \lambda^{-1} f^*(\eta^*) \ , \nonumber \\[-1.5ex]
& \\[-1.5ex]
& {\ds \frac{d f}{d \eta}(\eta) = \lambda^{\delta -1} \frac{d f^*}{d \eta^{*}}(\eta^*)} \ ,   
\quad
{\ds \frac{d^2 f}{d \eta^{2}}(\eta) = \lambda^{2\delta-1} \frac{d^2 f^*}{d \eta^{*2}}(\eta^*)} \nonumber \ .  
\end{align}
In particular, we are interested to compute the missing initial condition ${\frac{d^2 f}{d \eta^{2}}(0)}$.

We are now ready to present the method of solution in the form of an algorithm.

\noindent
{\bf The non-iterative algorithm.}\\
1. Input $p$, $\delta$, $\eta_{\infty}^*$, $d$.\\
2. Solve  (\ref{eq:IVP1}) in $[0, \eta_{\infty}^*]$.\\
3. Compute $\lambda$ by (\ref{eq:lambda1}).\\
4. Rescale the numerical solution according to (\ref{eq:rescale1}).\\

The above algorithm defines a non-ITM for the numerical solution of the class of problems characterized by (\ref{eq:condition}) and $a=b=c=0$.

\section{The extension due to Na}

We report here an important extension of the non-ITM due to Na.
For this extension we have to require the invariance of physical parameters.
The boundary conditions at $\eta = 0$ involve the parameters $a$, $b$ and $c$.
Those boundary conditions are left invariant by the scaling transformation on condition
that the involved parameters transform as
\begin{equation}\label{eq:ex:Na}
a^* = \lambda a \ , \quad b^* = \lambda^{1-\delta} b \ , \quad c^* = \lambda^{\delta} c \ .
\end{equation}

\subsection{Moving wall}

According to Ishak et al. \cite{Ishak:2007:BLM} the moving wall boundary conditions are given by
\begin{equation}\label{eq:BCs2}
f(0) = 0 \ , \qquad
\frac{df}{d\eta} (0) = b  \ , \qquad
{\displaystyle \frac{df}{d\eta}}(\eta) \rightarrow d \quad \mbox{as}
\quad \eta \rightarrow \infty \ ,
\end{equation}
where $d = 1-b$ and $ b $ is a non-dimensional parameter given by the ration of the wall to the flow velocities. 
The boundary conditions of the Blasius problem are recovered from (\ref{eq:BCs2}) by setting $b=0$.

\subsubsection{The non-ITM} 

The applicability of a non-ITM to the Blasius problem is a consequence of its invariance with respect to a scaling transformation; note that the asymptotic boundary condition is not invariant.
In order to apply the non-ITM we consider $ b $ as a parameter involved in the scaling invariance, i.e., we define the extended scaling group
\begin{equation}\label{eq:scalinv:wall}
f^* = \lambda f \ , \qquad \eta^* = \lambda^{-1} \eta \ , \qquad 
b^* = \lambda^{2} b \ .   
\end{equation}
By setting a value of $P^*$, we can integrate the Blasius equation (\ref{eq:Blasius}) in the star variables with initial conditions
\begin{equation}\label{eq:ICs2}
f^*(0) = 0 \ , \quad \frac{df^*}{d\eta^*}(0) = b^* \ , \quad \frac{d^2f^*}{d\eta^{*2}}(0) = \pm 1 \ ,
\end{equation}
in order to compute an approximation $ \frac{d f^*}{d \eta^{*}}(\eta_{\infty}) $ for $\frac{df^*}{d\eta^*}(\infty)$ and the corresponding value of $\lambda$ according to the equation %(\ref{eq:lambda:wall})
\begin{equation}\label{eq:lambda:wall}
\lambda = \left[ \frac{d f^*}{d \eta^{*}}(\eta_{\infty})+b^* \right]^{1/2} \ .   
\end{equation} 
Once the value of $\lambda$ is computed by equation (\ref{eq:lambda:wall}), then we can find the missed initial condition by the equation
\begin{equation}\label{eq:MIC1}
\frac{d^2f}{d\eta^{2}}(0) =  \lambda^{-3}\frac{d^2f^*}{d\eta^{*2}}(0) \ ,
\end{equation}

For the application of the method defined above, we remark that the plus (for $b<0.5$) or minus (when $b>0.5$) sign must be used for the second derivative in (\ref{eq:ICs2}).
Moreover, the computation of a value at infinity is unsuitable from a numerical viewpoint and therefore we use a truncated boundary $\eta^*_{\infty}$ instead of infinity.
\begin{figure}[!hbt]
	\centering
%\begin{small}
\psfrag{r}[][]{$b$} 
%\psfrag{df}[][]{$\frac{df}{d\eta}(0)$, $\frac{d^2f}{d\eta^{2}}(0)$} 
%\psfrag{ddff}[][]{} 
\psfrag{df-}[][r]{$\frac{df^*}{d\eta^*}(0)=-1$} 
\psfrag{df+}[l][l]{$\frac{df^*}{d\eta^*}(0)=1$} 
\psfrag{d2f}[][]{$\frac{d^2f}{d\eta^2}(0)$} 
\includegraphics[width=\tw]{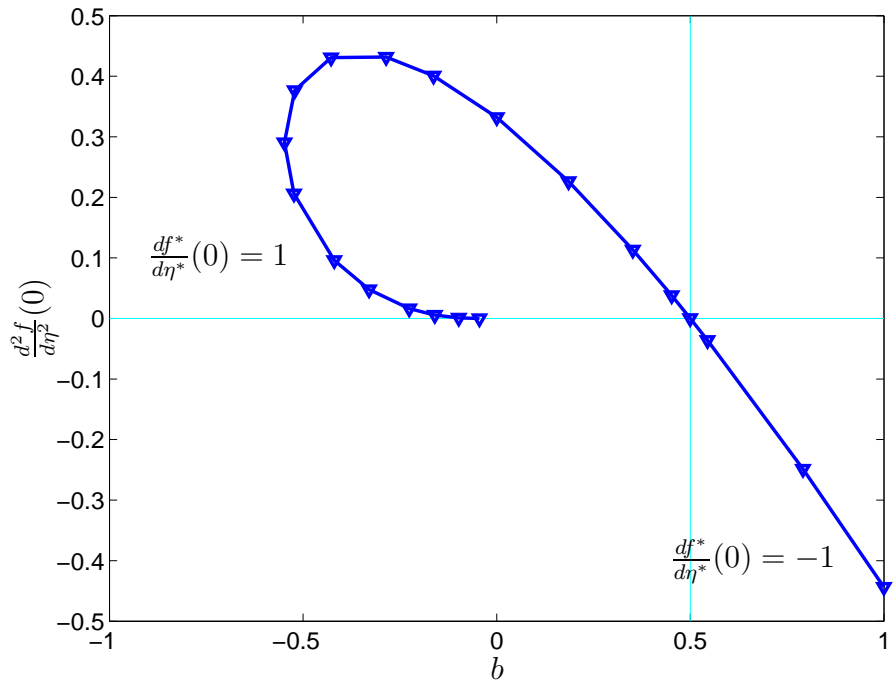}
%\end{small}
\caption{Blasius problem with moving wall boundary conditions. Non-ITM: missing initial condition versus $b$.}
	\label{fig:P-Mwall}
\end{figure}
In figure \ref{fig:P-Mwall} we plot $\frac{d^2f}{d\eta^2}(0)$ versus $b$.
From this figure we realize that our problem has an unique solution for $b\ge 0$, while two solutions exist for $b<0$.
The values of the skin friction coefficient are positive for $b<1/2$ and they become negative when $b>1/2$.
From a physical viewpoint, this means that a drag force is exerted by the flow on the plate in the first case, while in the second the force is of opposite type.
Figure \ref{fig:P-Mwall} also shows that the solutions exist until a critical, negative, value of the parameter, say $b_c$, is reached.
The boundary layer separate from the surface at $b=b_c$, and, therefore, for smaller values of $b$ the Navier-Stoker equations has to be solved because the hypotheses of boundary layer theory felt down.
We have separation for a positive value of the skin friction coefficient and not at the point where this coefficient vanish as in the classical boundary layer theory. 
The zero value of the skin friction coefficient when $b=1/2$ corresponds to equal velocity of the plate and the free stream and does not mark separation.

From the data in table \ref{tab:PMwall} we get $b_c \approx -0.548210$.   
This value is in good agreement with the value $b_c = -0.5483$ computed by Ishak et al. \cite{Ishak:2007:BLM} using an iterative method: the second order Keller's Box finite difference method.
%Table \ref{tab:PMwall} lists some numerical results.
\begin{table}[!htb]
%\begin{table}[p]
\renewcommand\arraystretch{1.1}
	\centering
		\begin{tabular}{cr@{.}lr@{.}lr@{.}lr@{.}l}
\hline \\[-3ex]
{${\displaystyle \frac{d^2f^*}{d{\eta^*}^2}(0)}$}
& \multicolumn{2}{c}%
{$ b^* $}
& \multicolumn{2}{c}%
{${\displaystyle \frac{df^*}{d\eta^*}(\infty)}$}
& \multicolumn{2}{c}%
{${\displaystyle \frac{d^2f}{d\eta^2}(0)}$}
& \multicolumn{2}{c}%
{$ b $} \\[2ex]
\hline
$1$ &
$-500$ & & $1$ & $55\cdot 10^{4}$ & $5$ & $46\cdot 10^{-7}$ & $-0$ & $033393$ \\
&$-100$ & & $2$ & $34\cdot 10^{3}$ & $9$ & $42\cdot 10^{-6}$ & $-0$ & $044591$ \\
&$-5$ &     & $36$ & $325698$ & $0$ & $005704$ & $-0$ & $159613$ \\
%&$-1$ & $5$  & $4$ & $368544$ & $0$ & $205830$ & $-0$ & $522913$ \\
%&$-1$ & $25$ & $3$ & $529165$ & $0$ & $290627$ & $-0$ & $548447$ \\
&$-1$ &      & $2$ & $917762$ & $0$ & $376537$ & $-0$ & $521441$ \\
%&$-0$ & $75$ & $2$ & $503099$ & $0$ & $430814$ & $-0$ & $427814$ \\
%&$-0$ & $5$  & $2$ & $250439$ & $0$ & $431797$ & $-0$ & $285643$ \\
&$0$ & & $2$ & $085393$ & $0$ & $332061$ & $0$ &  \\
%&$0$ & & $2$ & $085623$ & $0$ & $332006$ & $0$ & $$ \\
&$1$ & & $2$ & $440648$ & $0$ & $156689$ & $0$ & $290643$ \\
&$5$ & & $5$ & $771518$ & $0$ & $028287$ & $0$ & $464187$ \\
&$100$ & & $1$ & $00\cdot 10^{2}$ & $3$ & $53\cdot 10^{-4}$ & $0$ & $499557$ \\
&$500$ & & $5$ & $00\cdot 10^{2}$ & $3$ & $16\cdot 10^{-5}$ & $0$ & $499960$ \\
$-1$ &
$100$ & & $99$ & $822681$ & $-3$ & $54\cdot 10^{-4}$ & $0$ & $500444$ \\
&$10$ & & $9$ & $433763$ & $-0$ & $011673$ & $0$ & $514568$ \\
&$5$  & & $4$ & $182424$ & $-0$ & $035939$ & $0$ & $544519$ \\
&$2$  & & $0$ & $528464$ & $-0$ & $248722$ & $0$ & $790994$ \\
&$1$  & $719$ & $-4$ & $73\cdot 10{-5}$ & $-0$ & $443715$ & $1$ & $000027$ \\
\hline			
		\end{tabular}
	\caption{Moving wall boundary condition: non-ITM numerical results.}
	\label{tab:PMwall}
\end{table}
As mentioned before, the case $b=0$ is the Blasius problem (\ref{eq:Blasius}).
In this case our non-ITM becomes the original method defined by T\"opfer \cite{Topfer:1912:BAB}.
The obtained skin friction coefficient is in good agreement with the values available in literature, see for instance Fazio \cite{Fazio:1992:BPF} or Boyd \cite{Boyd:1999:BFC}. 
The values reported in the last line are related to the Sakiadis problem \cite{Sakiadis:1961:BLBa,Sakiadis:1961:BLBb} and were found by a few trial and miss attempts.
The obtained skin friction coefficient is in good agreement with the values reported in literature, e.g.: Sakiadis \cite{Sakiadis:1961:BLBa}, Ishak et al. \cite{Ishak:2007:BLM}, Cortell \cite{Cortell:2010:NCB}and Fazio \cite{Fazio:2015:ITM} with an iterative TM .

\begin{figure}[!hbt]
	\centering
%\begin{small}
\psfrag{e}[][]{$\eta$} 
%\psfrag{ddff}[][]{} 
\psfrag{df}[l][]{$\frac{df}{d\eta}$} 
\psfrag{ddf}[l][]{$\frac{d^2f}{d\eta^2}$} 
\psfrag{df*}[l][]{$\frac{df^*}{d\eta^*}$} 
\psfrag{ddf*}[l][]{$\frac{d^2f^*}{d{\eta^*}^2}$} 
\includegraphics[width=\tw,height=8cm]{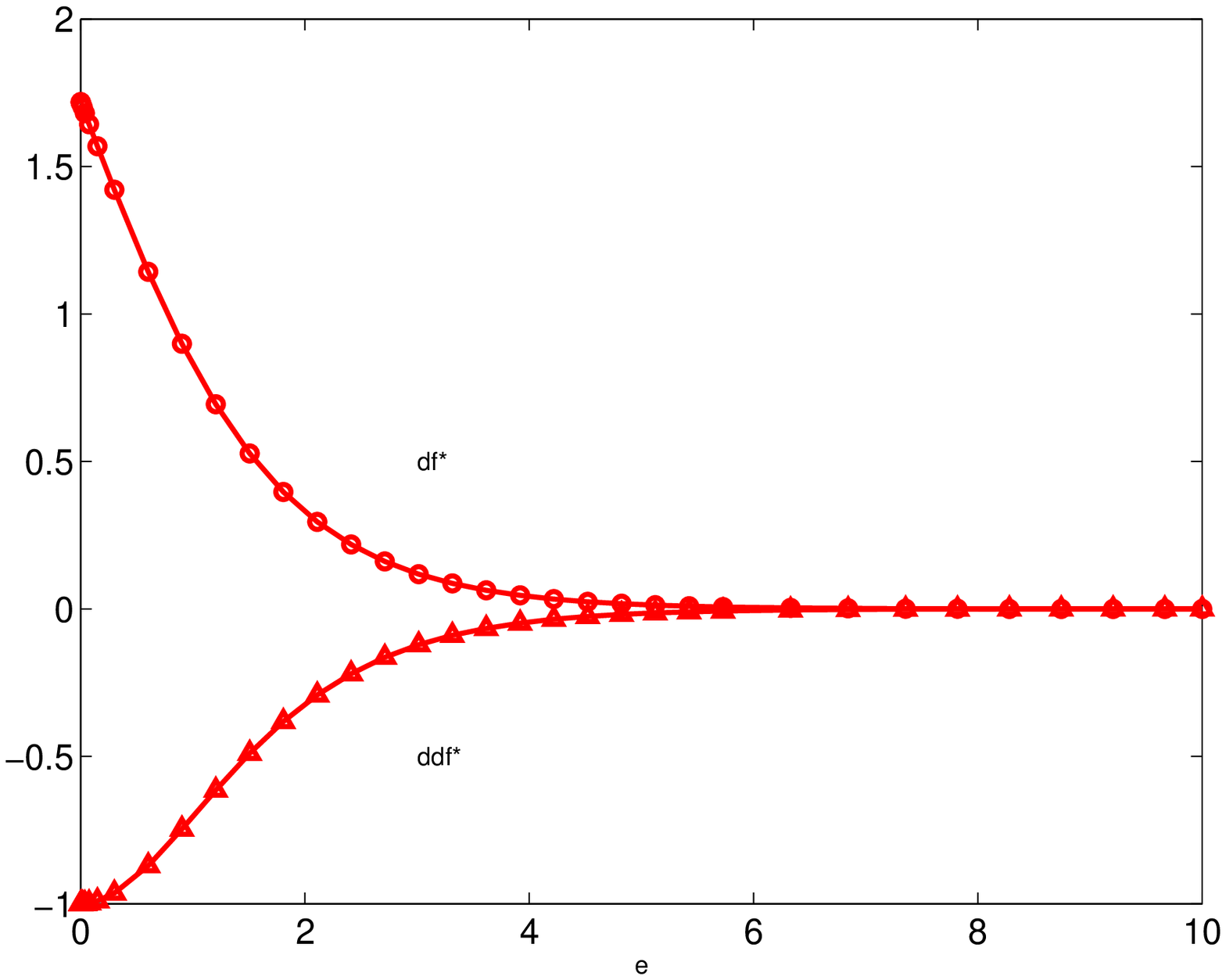} \\
\includegraphics[width=\tw,height=8cm]{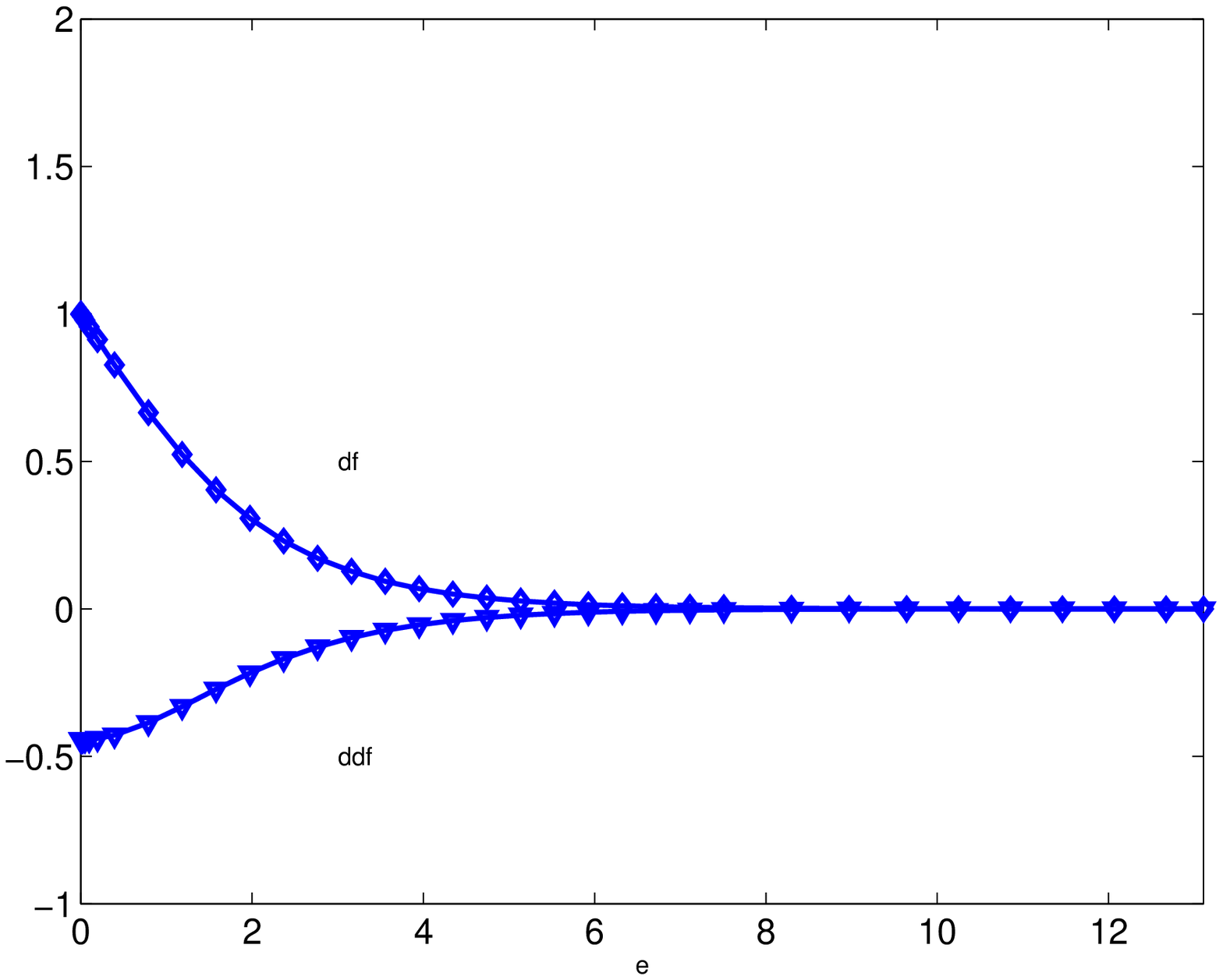}
%\end{small}
\caption{Numerical results of the non-ITM. Top frame: solution of the IVP; bottom frame: solution of the Sakiadis problem found after rescaling.}
	\label{fig:Sakiadis:fig}
\end{figure}

Figure \ref{fig:Sakiadis:fig} shows the solution of the Sakiadis problem, describing the behaviour of a boundary layer flow due to a moving flat surface immersed in an otherwise quiescent fluid, corresponding to $b=1$.
Actually, this is a case of practical interest if we are considering the plate as an idealization of an airplane wing.

\subsection{Slip flow condition}

We consider now the case of a rarefied flow where the no-slip condition at the wall, considered in the previous section, must be replaced by a slip-flow condition, see for instance Gad-el-Hak \cite{Gad-el-Hak:1999:FMM}.
For an isothermal wall, the slip condition can be defined as
\begin{equation}
v(y, 0) = \frac{2-\sigma}{\sigma} \ell \frac{\partial v}{\partial z} (y, 0) \ ,
\end{equation} 
where $ \ell $ is the mean free path, and $ \sigma $ is the tangential momentum accommodation coefficient.
Within a similarity transformation this slip boundary condition becomes
\begin{equation}
\frac{df}{d\eta} (0) = c \; \frac{d^2f}{d\eta^2} (0) \ ,
\end{equation}
where $ c $ is a non-dimensional parameter, that takes into account the behaviour at the surface, defined by
\begin{equation}
c = \frac{2-\sigma}{\sigma} Kn \; Re \; y^{1/2} \ ,
\end{equation}
where $ Kn $ and $ Re $ are the Knudsen and Reynolds numbers based on $ y $.

%Existence and Uniqueness for $ P > 0 $.
For the Blasius problem with slip condition we implemented the extended non-ITM.

\subsubsection{The non-ITM} 

In order to apply the non-ITM we consider $ c $ as a parameter involved in the scaling invariance, i.e., we defined the extended scaling group
\begin{equation}\label{eq:scalinv:slip1}
f^* = \lambda f \ , \qquad \eta^* = \lambda^{-1} \eta \ , \qquad 
c^* = \lambda^{-1} c \ .   
\end{equation}
Henceforth, $ \lambda $ is defined, once again, by equation (\ref{eq:lambda:Blasius}).

Sample numerical results are reported on table \ref{tab:P}.
\begin{table}[!htb]
%\begin{table}[p]
\renewcommand\arraystretch{1.5}
	\centering
		\begin{tabular}{r@{.}lr@{.}lr@{.}lr@{.}lr@{.}l}
\hline
\multicolumn{2}{c}%
{$ c^* $}
& \multicolumn{2}{c}%
{${\displaystyle \frac{df^*}{d\eta^*}(\infty)}$}
& \multicolumn{2}{c}%
{${\displaystyle \frac{df}{d\eta}(0)}$}
& \multicolumn{2}{c}%
{${\displaystyle \frac{d^2f}{d\eta^2}(0)}$}
& \multicolumn{2}{c}%
{$ c $} \\
\hline
0  &  & 2  & 085393 & 0 &        & 0 & 332061 & 0 &  \\
0  &1 & 2  & 090453 & 0 & 047836 & 0 & 330856 & 0 & 144584 \\
0  &5 & 2  & 191907 & 0 & 228112 & 0 & 308153 & 0 & 740255 \\
1  &  & 2  & 440648 & 0 & 409727 & 0 & 262266 & 1 & 562257 \\
%2  &  & 3  & 141625 & 0 & 636613 & 0 & 179584 &  &  \\
%3  &  & 3  & 968362 & 0 & 755979 & 0 & 126498 &  &  \\
%4  &  & 4  & 853857 & 0 & 824087 & 0 & 093512 &  &  \\
5  &  & 5  & 771518 & 0 & 866323 & 0 & 072122 &12 & 011992 \\
%6  &  & 6  & 708699 & 0 & 894361 & 0 & 057550 &  &  \\
%7  &  & 7  & 658853 & 0 & 913975 & 0 & 047180 &  &  \\
%8  &  & 8  & 618107 & 0 & 928278 & 0 & 039526 &  &  \\
%9  &  & 9  & 583812 & 0 & 939084 & 0 & 033705 &  &  \\
10 &  & 10 & 554805 & 0 & 947436 & 0 & 029162 &32 & 488159 \\
15 &  & 15 & 455238 & 0 & 970545 & 0 & 016458 & 3 & 815517 \\
20 &  & 20 & 394883 & 0 & 980638 & 0 & 010857 &90 & 321389 \\
25 &  & 25 & 353618 & 0 & 986053 & 0 & 007833 &125& 880941 \\
\hline			
		\end{tabular}
	\caption{Slip boundary condition: non-iterative numerical results.}
	\label{tab:P}
\end{table}
As mentioned before, the case $c=0$ is the classical Blasius problem (\ref{eq:Blasius}).
It is possible to compare the results listed in the last two columns of table \ref{tab:P} with similar results, obtained via a shooting method, shown in figure 1 of the proceedings report by Martin and Boyd \cite{Martin:2001:BBL}.
It is clear that our non-ITM would be faster and easier to implement than any iterative algorithm.

\begin{figure}[!hbt]
\psfrag{e}[][]{$ \eta^*, \eta $} 
\psfrag{}[][]{$  $} 
\psfrag{f2star}[][]{$ \frac{df^*}{d\eta^*} $} 
\psfrag{f2}[][]{$ \frac{df}{d\eta} $} 
\psfrag{f3star}[][]{$ \frac{d^2f^*}{d\eta^{*2}} $} 
\psfrag{f3}[][]{$ \frac{d^2f}{d\eta^2} $} 
\centerline{
\includegraphics[width=\tw,height=8cm]{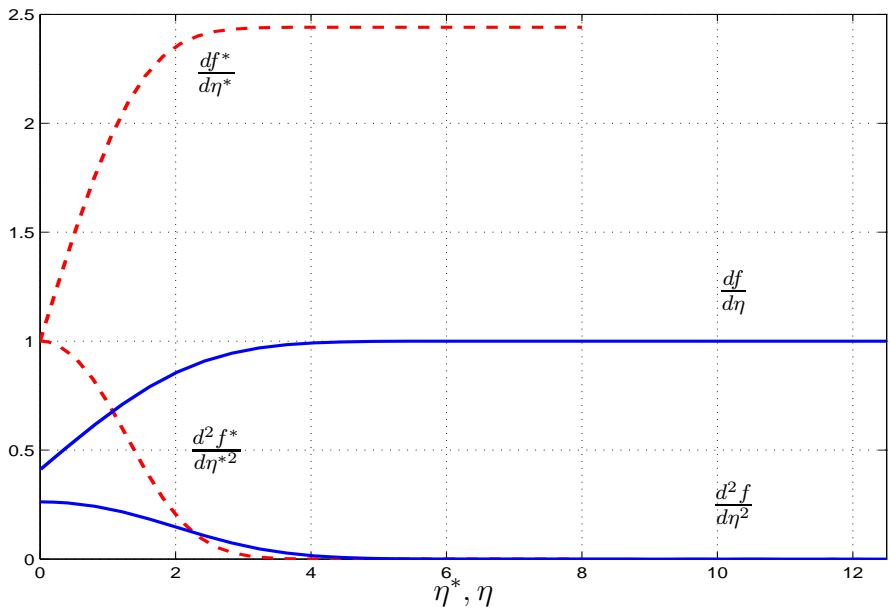}
}
\caption{Blasius problem with slip condition. Numerical solution by a non-ITM with $ c^* = 1 $ and $ c = 1.562257 $.}
	\label{fig:Blasius-Slip}
\end{figure}

As far as the non-ITM is concerned, we set a value of $ c^* $ and get the numerical solution of the problem for a different value of $c$.
As an example, figure \ref{fig:Blasius-Slip} shows a sample numerical integration for $ c = 1.562257 $ obtained by fixing $c^* = 1$.
Note that the solution of the Blasius problem with slip boundary condition was computed by rescaling.
If we need the solution for a specific value of $c$, then we can apply interpolation techniques to the results of table \ref{tab:P}. 

\begin{figure}[!hbt]
	\centering
\psfrag{P}[][]{$c$} 
%\psfrag{df}[][]{$\frac{df}{d\eta}(0)$, $\frac{d^2f}{d\eta^{2}}(0)$} 
\psfrag{ddff}[][]{} 
\psfrag{df}[][]{$\frac{df}{d\eta}(0)$} 
\psfrag{d2f}[][]{$\frac{d^2f}{d\eta^2}(0)$} 
\includegraphics[width=\tw,height=8cm]{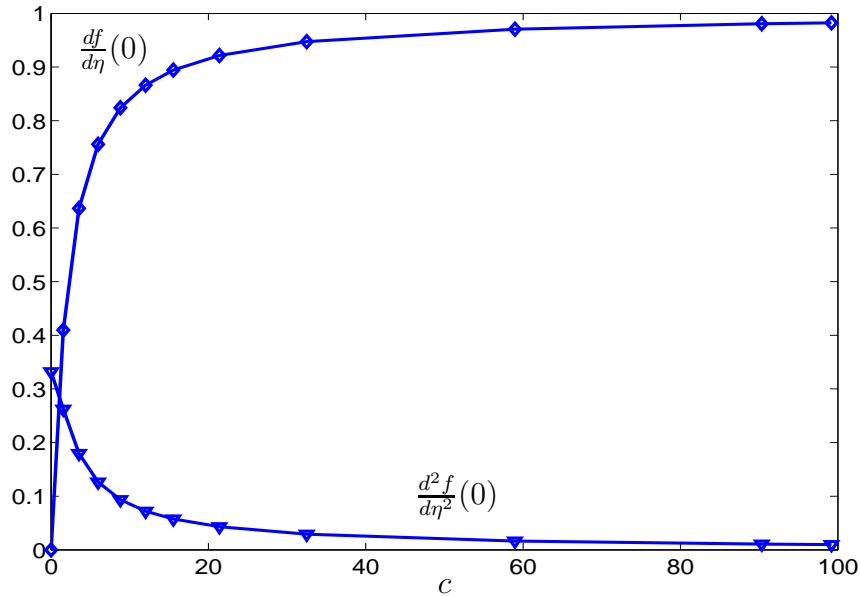}
%\end{small}
\caption{Blasius problem with slip boundary condition. Non-ITM: behaviour of $\frac{df}{d\eta}(0)$ and $\frac{d^2f}{d\eta^2}(0)$ with respect to $c$.}
	\label{fig:P-Slip}
\end{figure}

In figure \ref{fig:P-Slip} we plot $\frac{df}{d\eta}(0)$ and $\frac{d^2f}{d\eta^2}(0)$ versus $c$.
It is easily seen that as $c$ goes to infinity then $\frac{df}{d\eta}(0)$ goes to one while $\frac{d^2f}{d\eta^2}(0)$ tends to zero.

\section{Surface gasification boundary condition}

In this section we consider a problem that does not belong to the framework used so far.
This is the surface gasification flow that, with and without subsequent gas-phase flame-sheet burning, was first formulated and analysed by Emmons \cite{Emmons:1956:FCL}.
For such a flow we have to consider the variant of the celebrated Blasius problem
\begin{align}\label{eq:SGas}
& {\displaystyle \frac{d^3 f}{d \eta^3}} + f
{\displaystyle \frac{d^{2}f}{d\eta^2}}
= 0 \nonumber \\[-1ex]
&\\[-1ex]
& f(0) = -s \frac{d^2f}{d\eta^2} (0) \ , \qquad
\frac{df}{d\eta} (0) = 0 \ , \qquad
{\displaystyle \frac{df}{d\eta}}(\eta) \rightarrow 1 \quad \mbox{as}
\quad \eta \rightarrow \infty \ , \nonumber
\end{align}
where $ s $ is the classical Spalding heat transfer number \cite{Spalding:1961:MTT}. 
This transfer number for slow vaporization belongs to the interval $[0, 0.1]$ and varies from $s = O(1)$ to $s \approx 20$ for strong burning.
This problem has been studied recently by Lu and Law \cite{Lu:2014:ISB}.
These authors define an iterative method that has been shown to produce more accurate numerical results than the classical approximate solutions.

\subsection{The non-ITM} 

%In order to apply the non-ITM we consider $ P $ as a parameter involved in the scaling invariance, i.e.,
In the present case we consider the extended scaling group
\begin{equation}\label{eq:scalinv:gass1}
f^* = \lambda f \ , \qquad \eta^* = \lambda^{-1} \eta \ , \qquad 
s^* = \lambda^{-2} s \ .   
\end{equation}
Let us notice that the governing differential equation and the two boundary conditions at $\eta = 0$ in (\ref{eq:SGas}) are left invariant under (\ref{eq:scalinv:gass1}) and, on the contrary, the asymptotic boundary condition is not invariant.
By setting a value of $s^*$, we can integrate the Blasius governing differential equation in (\ref{eq:SGas}) in the star variables on $[0, \eta^*_\infty]$ with initial conditions
\begin{equation}\label{eq:ICs:gass}
f^*(0) = -s^* \ , \quad \frac{df^*}{d\eta^*}(0) = 0 \ , \quad \frac{d^2f^*}{d\eta^{*2}}(0) = 1 \ ,
\end{equation}
in order to compute $\frac{df^*}{d\eta^*}(\eta^*_\infty) \approx \frac{df^*}{d\eta^*}(\infty)$.
Here $\eta^*_\infty$ is a suitable truncated boundary.
The value of $ \lambda $ can be found by
\begin{equation}\label{eq:lambda4}
\lambda = \left[ \frac{d f^*}{d \eta^{*}}(\eta_\infty^*) \right]^{1/2} \ .   
\end{equation} 
After using (\ref{eq:lambda4}) to get the value of $\lambda$, we can apply the scaling invariance to obtain the missing initial conditions
\begin{equation}\label{eq:MICs:gas}
f(0) = \lambda^{-2} s^* \ , \quad \frac{d^2f}{d\eta^{2}}(0) = \lambda^{-3} \ .
\end{equation}

For the reader convenience, in table \ref{tab:P-Gas} we list sample numerical results.
\begin{table}[!htb]
%\begin{table}[p]
\renewcommand\arraystretch{1.1}
	\centering
		\begin{tabular}{r@{.}lr@{.}lr@{.}lr@{.}lr@{.}l}
\hline \\[-3ex]
\multicolumn{2}{c}%
{$ s^* $}
& \multicolumn{2}{c}%
{${\displaystyle \frac{df^*}{d\eta^*}(\infty)}$}
& \multicolumn{2}{c}%
{${-f(0)}$}
& \multicolumn{2}{c}%
{${\displaystyle \frac{d^2f}{d\eta^2}(0)}$}
& \multicolumn{2}{c}%
{$ s $} \\[2ex]
\hline
0  &    & 1  & 655301 & 0 &        & 0 & 469553 & 0 &         \\
%0  & 1  & 1  & 793644 & 0 & 074668 & 0 & 416289 & 0 & 179364  \\
0  & 25 & 2  & 025902 & 0 & 175643 & 0 & 346795 & 0 & 506476  \\
0  & 5  & 2  & 485809 & 0 & 317129 & 0 & 255152 & 1 & 242904  \\
0  & 75 & 3  & 048481 & 0 & 429556 & 0 & 187877 & 2 & 286361  \\
1  &    & 3  & 726397 & 0 & 518031 & 0 & 139016 & 3 & 726397  \\
1  & 25 & 4  & 528469 & 0 & 587401 & 0 & 103770 & 5 & 660586  \\
1  & 5  & 5  & 469166 & 0 & 641403 & 0 & 078184 & 8 & 203749  \\
1  & 75 & 6  & 548781 & 0 & 683845 & 0 & 059670 &11 & 460366  \\
2  &    & 7  & 779561 & 0 & 717055 & 0 & 046086 &15 & 559122  \\
%2  & 2  & 8  & 863956 & 0 & 738939 & 0 & 037893 &19 & 500704  \\
\hline			
		\end{tabular}
	\caption{Surface gasification boundary condition: non-ITM results.}
	\label{tab:P-Gas}
\end{table}
The case $s^*=s=0$ is, again, the Blasius problem (\ref{eq:SGas}).
In this case our non-ITM reduces to the original method defined by T\"opfer \cite{Topfer:1912:BAB}.
The obtained skin friction coefficient is in good agreement with the values available in literature, see for instance the value $0.46599988361$ computed by Fazio \cite{Fazio:1992:BPF}. 
On the other hand, our value is different from the value $0.490$ obtained by a $3-2$ iteration solution of Lu and Law \cite{Lu:2014:ISB}.
For the numerical results reported here, depending on the behaviour of the numerical solution, we have used $\eta^*_\infty = 5$ or $\eta^*_\infty = 10$.

\begin{figure}[!hbt]
	\centering
%\begin{small}
\psfrag{e}[][]{$\eta^*$, $\eta$} 
%\psfrag{}[][]{$  $} 
\psfrag{f2star}[][]{$\frac{df^*}{d\eta^*}$} 
\psfrag{f2}[][]{$\frac{df}{d\eta}$} 
\psfrag{f3star}[][]{$\frac{d^2f^*}{d\eta^{*2}}$} 
\psfrag{f3}[][]{$\frac{d^2f}{d\eta^2}$} 
\includegraphics[width=\tw]{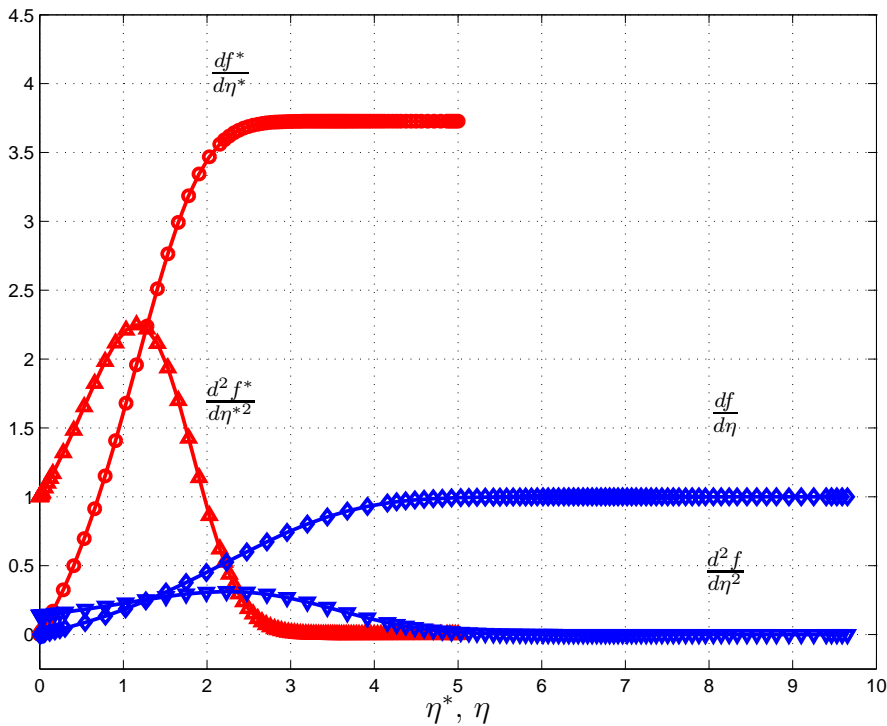}
%\end{small}
\caption{Surface gasification boundary conditions with $s^* = 1$. Numerical solution by the non-ITM.}
	\label{fig:Blasius-gas}
\end{figure}

Figure \ref{fig:Blasius-gas} shows a sample numerical integration for $ s^* = 1 $ that is transformed under (\ref{eq:scalinv:gass1}) to $s \approx 3.726397$.
We notice that the solution of the Blasius problem with surface gasification boundary condition is computed by rescaling.
Moreover, by rescaling we get $\eta^*_{\infty} < \eta_{\infty}$. 
%\begin{figure}[p]

\begin{figure}[!hbt]
	\centering
%\begin{small}
\psfrag{P}[][]{$s$} 
\psfrag{fdf}[][]{$f(0)$, $\frac{d^2f}{d\eta^2}(0)$} 
\psfrag{f}[][]{$f(0)$} 
\psfrag{d2f}[][]{$\frac{d^2f}{d\eta^2}(0)$} 
\includegraphics[width=\tw]{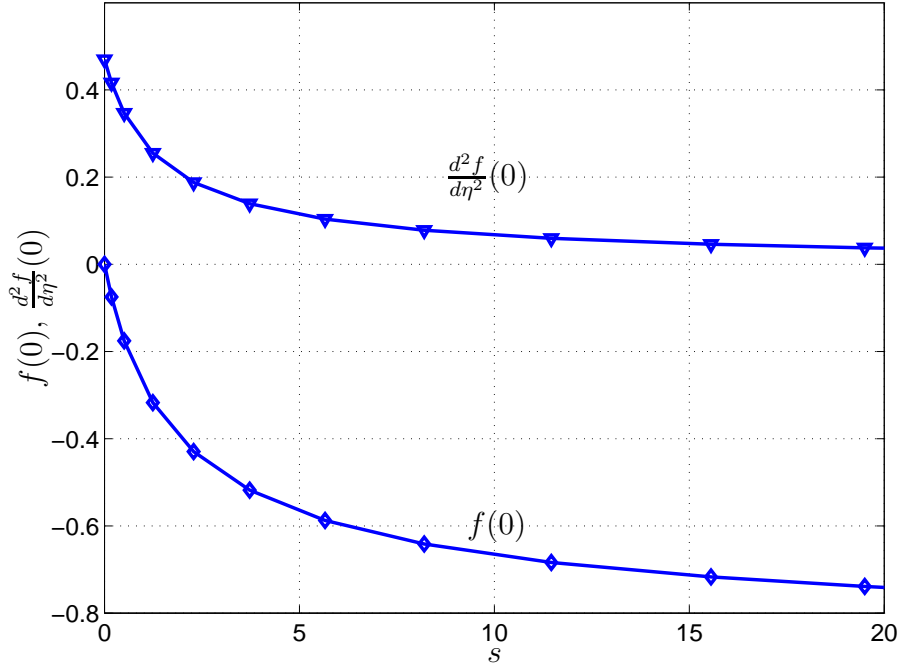}
%\end{small}
\caption{Blasius problem with surface gasification boundary conditions. Non-ITM: behaviour of $f(0)$ and $\frac{d^2f}{d\eta^2}(0)$ with respect to $s$.}
	\label{fig:P-gas}
\end{figure}

In figure \ref{fig:P-gas} we plot $f(0)$ and $\frac{d^2f}{d\eta^2}(0)$ versus $s$.
We know, from the literature, that as $s$ goes to infinity then $f(0)$ goes to $-0.876$. 
Moreover, as it is easily seen, as $s$ goes to infinity then $\frac{d^2f}{d\eta^2}(0)$ goes to zero. 

\section{Concluding remarks}

The main contribution of this paper is the extension of the non-ITM proposed by T\"opfer \cite{Topfer:1912:BAB} for the numerical solution of the celebrated Blasius problem \cite{Blasius:1908:GFK} to classes of problems depending on a parameter.
By requiring the invariance of the involved parameter we are able to solve the given BVP non-iteratively but for a different value of the parameter.
This kind of extension was considered first by Na \cite{Na:1970:IVM}, see also NA \cite[Chapters 8-9]{Na:1979:CME}.
Here we defined a non-ITM for Blasius equation with moving wall, surface gasification or slip boundary conditions.

Finally, the reader should be advised that non-ITM cannot be applied to all problems of boundary layer theory.
As an example, let us consider the Falkner-Skan model \cite{Falkner:1931:SAS}
\begin{align}\label{eq:abf}
& {\displaystyle \frac{d^{3}f}{d\eta^3}} + f 
{\displaystyle \frac{d^{2}f}{d\eta^2}} + P \; \left[1 - \left({\displaystyle
\frac{df}{d\eta}}\right)^2 \right] = 0 \ , \nonumber \\[-1.2ex]
& \\[-1.2ex]
& f(0) = {\displaystyle \frac{df}{d\eta}}(0) = 0 \ , \qquad
{\displaystyle \frac{df}{d\eta}}(\eta) \rightarrow 1 \quad \mbox{as}
\quad \eta \rightarrow \infty \ , \nonumber
\end{align}
where $f$ and $\eta$ are similarity variables and $P$ is a parameter related to the functional form of the fluid mainstream velocity.
If we test the invariance of the governing differential equation in (\ref{eq:abf}) under the extended scaling group
\begin{equation}\label{eq:sgext}
\eta^* = \lambda^{\alpha_1} \eta \ , \qquad f^* = \lambda^{\alpha_2} f \ , \qquad 
P^* = \lambda^{\alpha_3} P \ ,   
\end{equation}
where $\lambda$ is, again, the group parameter and ${\alpha_j}$, for $j = 1, 2, 3$, are constant to be determined, then we get three invariant conditions
\begin{equation}\label{eq:icond}
\alpha_2 -3 \alpha_1 = 2 (\alpha_2 - \alpha_1) = \alpha_3 = \alpha_3 + 2 (\alpha_2 - \alpha_1) \ .
\end{equation}
Now, it is a simple matter to show that the linear system defined by (\ref{eq:icond}) has the unique solution $\alpha_1 = \alpha_2 = \alpha_3 = 0$.
%the second equation in (\ref{eq:icond}) means $\alpha_3 = 2 (\alpha_2 - \alpha_1)$ and this replaced in the third equation imply $\alpha_3 = 2 \alpha_3$ that is $\alpha_3 = 0$. Therefore, we end up with the 
However, as a last word on this topic we can mention an iterative extension of our transformation method that has been developed in \cite{Fazio:1994:FSE,Fazio:1996:NAN} and successfully applied to the Falkner-Skan model \cite{Fazio:1994:FSE,Fazio:2013:BPF}.

%% Example of a theorem:
%\begin{Theorem}
%Text text text
%\end{Theorem}

%% Example of a proof:
%\begin{proof}[Proof of Theorem 1]
%Text text text
%\end{proof}

%%%%%%%%%%%%%%%%%%%%%%%%%%%%%%%%%%%%%%%%%%

%\acknowledgments{Acknowledgments}

%Main text.

%%%%%%%%%%%%%%%%%%%%%%%%%%%%%%%%%%%%%%%%%%

%\authorcontributions{Author Contributions}

%Required if more than one author. Authorship must include and be strictly limited to researchers who have substantially contributed to the reported work. Please carefully review our criteria regarding the Qualification for Authorship: \web /instructions.

%%%%%%%%%%%%%%%%%%%%%%%%%%%%%%%%%%%%%%%%%%

%=================================================================
% Back Matter (References and Notes)
%----------------------------------------------------------%=================================================================
% References:  Variant B
%=================================================================
% Use the following option to include external BibTeX files:
%\bibliography{lite}

\begin{thebibliography}{10}

\bibitem{Abbasbandy:2007:NSB}
S.~Abbasbandy.
\newblock A numerical solution of {B}lasius equation by {A}domian's
  decomposition method and comparison with homotopy perturbation method.
\newblock {\em Caos, Solitons \& Fractals}, 31:257--260, 2007.

\bibitem{Bairstow:1925:SF}
L.~Bairstow.
\newblock Skin friction.
\newblock {\em J. Roy. Aero. Soc.}, 29:3--23, 1925.

\bibitem{Cortell:2010:NCB}
R.~Cortell Bataller.
\newblock Numerical comparisons of {B}lasius and {S}akiadis flows.
\newblock {\em Matematika}, 26:187--196, 2010.

\bibitem{Blasius:1908:GFK}
H.~Blasius.
\newblock Grenzschichten in {F}l\"{u}ssigkeiten mit kleiner {R}eibung.
\newblock {\em Z. Math. Phys.}, 56:1--37, 1908.

\bibitem{Boyd:1999:BFC}
J.~P. Boyd.
\newblock The {B}lasius function in the complex plane.
\newblock {\em Exp. Math.}, 8:381--394, 1999.

\bibitem{Boyd:2008:BFC}
J.~P. Boyd.
\newblock The {B}lasius function: computation before computers, the value of
  tricks, undergradute projects, and open research problems.
\newblock {\em SIAM Rev.}, 50:791--804, 2008.

\bibitem{Collatz}
L.~Collatz.
\newblock {\em The Numerical Treatment of Differential Equations}.
\newblock Springer, Berlin, 3rd edition, 1960.

\bibitem{Gad-el-Hak:1999:FMM}
M.~Gad el~Hak.
\newblock The fluid mechanics of microdevices --- the {F}reeman scholar
  lecture.
\newblock {\em J. Fluids Eng.}, 121:5--33, 1999.

\bibitem{Emmons:1956:FCL}
H.~W. Emmons.
\newblock The film combustion of liquid fluid.
\newblock {\em ZAMM - J. Appl. Math. Mech.}, 36:60--71, 1956.

\bibitem{Falkner:1936:MNS}
V.~M. Falkner.
\newblock A method of numerical solution of differential equations.
\newblock {\em Philos. Mag.}, 21:624--640, 1936.

\bibitem{Falkner:1931:SAS}
V.~M. Falkner and S.~W. Skan.
\newblock Some approximate solutions of the boundary layer equations.
\newblock {\em Philos. Mag.}, 12:865--896, 1931.

\bibitem{Fazio:1992:BPF}
R.~Fazio.
\newblock The {B}lasius problem formulated as a free boundary value problem.
\newblock {\em Acta Mech.}, 95:1--7, 1992.

\bibitem{Fazio:1994:FSE}
R.~Fazio.
\newblock The {F}alkner-{S}kan equation: numerical solutions within group
  invariance theory.
\newblock {\em Calcolo}, 31:115--124, 1994.

\bibitem{Fazio:1996:NAN}
R.~Fazio.
\newblock A novel approach to the numerical solution of boundary value problems
  on infinite intervals.
\newblock {\em SIAM J. Numer. Anal.}, 33:1473--1483, 1996.

\bibitem{Fazio:2002:SFB}
R.~Fazio.
\newblock A survey on free boundary identification of the truncated boundary in
  numerical {BVP}s on infinite intervals.
\newblock {\em J. Comput. Appl. Math.}, 140:331--344, 2002.

\bibitem{Fazio:2013:BPF}
R.~Fazio.
\newblock {B}lasius problem and {F}alkner-{S}kan model: {T}{\"o}pfer's
  algorithm and its extension.
\newblock {\em Comput. \& Fluids}, 73:202--209, 2013.

\bibitem{Fazio:2015:ITM}
R.~Fazio.
\newblock The iterative transformation method for the {S}akiadis problem.
\newblock {\em Comput. \& Fluids}, 106:196--200, 2015.

\bibitem{Fox}
L.~Fox.
\newblock {\em Numerical Solution of Two-point Boundary Value Problems in
  Ordinary Differential Equations}.
\newblock Clarendon Press, Oxford, 1957.

\bibitem{Goldstein:1930:CSS}
S.~Goldstein.
\newblock Concerning some solutions of the boundary layer equations in
  hydro-dynamics.
\newblock {\em Proc. Camb. Philos. Soc.}, 26:1--30, 1930.

\bibitem{Horwarth:1938:SLB}
L.~Horwarth.
\newblock On the solution of the laminar boundary layer equations.
\newblock {\em Proc. Roy. Soc. London A}, 164:547--579, 1938.

\bibitem{Ishak:2007:BLM}
A.~Ishak, R.~Nazar, and I.~Pop.
\newblock Boundary layer on a moving wall with suction and injection.
\newblock {\em Chin. Phys. Lett.}, 24:2274--2276, 2007.

\bibitem{Lu:2014:ISB}
Z.~Lu and C.~K. Law.
\newblock An iterative solution of the {B}lasius flow with surface
  gasification.
\newblock {\em Int. J. Heat and Mass Transfer}, 69:223--229, 2014.

\bibitem{Martin:2001:BBL}
M.~J. {Martin} and I.~D. {Boyd}.
\newblock {Blasius boundary layer solution with slip flow conditions}.
\newblock In {\em Rarefied Gas Dynamics: 22nd International Symposium}, volume
  585 of {\em American Institute of Physics Conference Proceedings}, pages
  518--523, 2001, DOI: 10.1063/1.1407604.

\bibitem{Na:1970:IVM}
T.~Y. Na.
\newblock An initial value method for the solution of a class of nonlinear
  equations in fluid mechanics.
\newblock {\em J. Basic Engrg. Trans. ASME}, 92:503--509, 1970.

\bibitem{Na:1979:CME}
T.~Y. Na.
\newblock {\em Computational Methods in Engineering Boundary Value Problems}.
\newblock Academic Press, New York, 1979.

\bibitem{Prandtl:1904:UFK}
L.~Prandtl.
\newblock {\"U}ber {F}l{\"u}ssigkeiten mit kleiner {R}eibung.
\newblock In {\em Proceedings Third Internernatinal Math. Congress}, pages
  484--494, 1904.
\newblock Engl. transl. in {NACA} Tech. Memo. 452.

\bibitem{Rubel:1955:EET}
L.~A. Rubel.
\newblock An estimation of the error due to the truncated boundary in the
  numerical solution of the {B}lasius equation.
\newblock {\em Quart. Appl. Math.}, 13:203--206, 1955.

\bibitem{Sakiadis:1961:BLBa}
B.~C. Sakiadis.
\newblock Boundary-layer behaviour on continuous solid surfaces: I.
  {B}oundary-layer equations for two-dimensional and axisymmetric flow.
\newblock {\em AIChE J.}, 7:26--28, 1961.

\bibitem{Sakiadis:1961:BLBb}
B.~C. Sakiadis.
\newblock Boundary-layer behaviour on continuous solid surfaces: {II}. {T}he
  boundary layer on a continuous flat surface.
\newblock {\em AIChE J.}, 7:221--225, 1961.

\bibitem{Schlichting:1979:BLT}
H.~Schlichting.
\newblock {\em Boundary-Layer Theory}.
\newblock McGraw-Hill, New York, 1979.

\bibitem{Spalding:1961:MTT}
D.~B. Spalding.
\newblock Mass transfer through laminar boundary layers - 1. {T}he velocity
  boundary layer.
\newblock {\em Int. J. Heat Mass Transfer}, 2:15--32, 1961.

\bibitem{Tajvidi:2008:MRL}
M.~Tajvidi, M.~Razzaghi, and M.~Dehghan.
\newblock Modified rational {L}egendre approach to laminar viscous flow over a
  semi-infinite flat plate.
\newblock {\em Caos, Solitons \& Fractals}, 35:59--66, 2008.

\bibitem{Topfer:1912:BAB}
K.~T{\"o}pfer.
\newblock Bemerkung zu dem {A}ufsatz von {H}. {B}lasius: {G}renzschichten in
  {F}l{\"u}ssigkeiten mit kleiner {R}eibung.
\newblock {\em Z. Math. Phys.}, 60:397--398, 1912.

\bibitem{Wazwaz:2007:VIM}
A.-M. Wazwaz.
\newblock The variational iteration method for solving two forms of {B}lasius
  equation on a half-infinite domain.
\newblock {\em Appl. Math. Comput.}, 188:485--491, 2007.

\bibitem{Weyl:1942:DES}
H.~Weyl.
\newblock On the differential equation of the simplest boundary-layer problems.
\newblock {\em Ann. Math.}, 43:381--407, 1942.

\end{thebibliography}

%%%%%%%%%%%%%%%%%%%%%%%%%%%%%%%%%%%%%%%%%%

%\abbreviations{Abbreviations/Nomenclature}
%
%Main text.

%%%%%%%%%%%%%%%%%%%%%%%%%%%%%%%%%%%%%%%%%%

%\appendix
%\section{Appendix Title}
%
%Main text.

\end{document}